# NUMERICAL SOLUTION OF OBSTACLE SCATTERING PROBLEMS

**Alexander G. Ramm**

*Department of Mathematics, Kansas State University*
*Manhattan, Kansas 66506-2602, USA, email: ramm@math.ksu.edu*

**Semion Gutman**

*Department of Mathematics, University of Oklahoma*
*Norman, OK 73019, USA, email: sgutman@ou.edu*

**Abstract**

Some novel numerical approaches to solving direct and inverse obstacle scattering problems (IOSP) are presented. Scattering by finite obstacles and by periodic structures is considered. The emphasis for solving direct scattering problem is on the Modified Rayleigh Conjecture (MRC) method, recently introduced and tested by the authors. This method is used numerically in scattering by finite obstacles and by periodic structures. Numerical results it produces are very encouraging. The support function method (SFM) for solving the IOSP is described and tested in some examples. Analysis of the various versions of linear sampling methods for solving IOSP is given and the limitations of these methods are described.

## 1. Introduction

In this paper we review our recent work on Direct and Inverse Obstacle Scattering Problems (IOSP), which is centered on the Rayleigh Conjecture (RC). These results show that numerical solution of the obstacle scattering problem based on the MRC (Modified Rayleigh Conjecture) method is a competitive alternative to the BIEM (boundary integral equations method). This approach has numerical advantages which may be especially important in three-dimensional scattering problems, and non-smooth domains, for example, in domains whose boundaries contain corners in 2-D case, and vertices and edges in 3-D case. Recently there was an increased interest in scattering by periodic structures. We discuss an adaptation of the MRC method to such scattering problems after a short introduction into its theory. We also present a novel method (Support Function Method) for the solution of the Inverse Obstacle Scattering Problem (IOSP), and compare its performance to the Linear Sampling Method.

We start with the formulation of the obstacle scattering problem. In this paper we usually consider the 2-D setting, and the Dirichlet boundary condition, but the discussed methods can also be used for 3-D problems, and the Neumann and Robin boundary conditions as well.

Let an obstacle be a bounded domain $D \subset \mathbb{R}^2$ with a Lipschitz boundary $\Gamma$. Fix a frequency $k > 0$ and denote the exterior domain by $D' = \mathbb{R}^2 \setminus \overline{D}$.

A solution $v(x)$ of the Helmholtz equation

$$\Delta v + k^2 v = 0, \quad x \in D', \quad (1.1)$$

is called outgoing if it satisfies the radiation condition

$$\lim_{r \to \infty} \int_{|x|=r} \left| \frac{\partial v}{\partial |x|} - ikv \right|^2 ds = 0. \quad (1.2)$$

The **Exterior Dirichlet Problem** consists of finding an outgoing solution of the Helmholtz equation (1.1) which satisfies the boundary condition

$$v = f, \quad x \in \Gamma, \qquad (1.3)$$

where $f$ is a continuous function, see [23] for the existence and uniqueness results for this problem, and [43], [45] for such results in the case of very rough domains. In [34] uniqueness of the solution to IOSP is proved for domains with finite perimeter, and the separability of $L^2(D)$ is used in the proof of the uniqueness results in place of the discreteness of the spectrum of Neumann Laplacian, thus showing the applicability of the Schiffer's beautiful idea to the IOSP uniqueness proof when the boundary condition is the Neumann one. In [41], [36], [37], [38], [39] and [40] various problems related to IOSP are studied: the dependence of the scattering amplitude on the boundary $S$ of the obstacle, and the stability of the obstacle towards small perturbation of the scattering amplitude.

A particular case of the above problem is the **Direct Acoustic Obstacle Scattering Problem**. Let $\alpha \in S^1$, and the incident field be

$$u^i(x) = e^{ikx \cdot \alpha}. \qquad (1.4)$$

The problem is to find the total field

$$u(x, k) = u^i + v, \quad x \in D' \qquad (1.5)$$

such that

$$u = 0, \quad x \in \Gamma, \qquad (1.6)$$

and the scattered field $v := u^s$ satisfies (1.1)-(1.2).

It is known (see e.g. [23]), that every outgoing solution $v(x)$, $x \in D'$ has an asymptotic representation

$$v(x) = \frac{e^{ik|x|}}{\sqrt{|x|}} \left\{ A(\alpha') + O\left(\frac{1}{|x|}\right) \right\}, \quad |x| \to \infty ,$$

(1.7)

where $\alpha' = x/|x|$, $\alpha' \in S^1$. The function $A(\alpha') := A_v(\alpha')$ is called the Far Field Pattern of $u$.

For the Direct Acoustic Obstacle Scattering Problem this representation takes the form

$$v(x) = \frac{e^{ik|x|}}{\sqrt{|x|}} \left\{ A(\alpha', \alpha) + O\left(\frac{1}{|x|}\right) \right\}, \quad |x| \to \infty ,$$

(1.8)

where the uniquely defined function $A(\alpha', \alpha)$ is called the Scattering Amplitude of the Obstacle Scattering Problem.

Let $J_l(t)$ and $N_l(t)$ be the Bessel and Neumann functions of the integer order $l$. The first Hankel function of order $l$ is defined by $H_l^{(1)} = J_l + iN_l$. Suppose that the circle $B_R = \{x \in \mathbb{R}^2 : |x| \leq R\}$ contains $D$. Then, in the region $|x| > R$, the outgoing solution of the Exterior Dirichlet Problem (1.1)-(1.3) has a unique representation

$$v(x) = \sum_{l=-\infty}^{\infty} a_l H_l^{(1)}(k|x|) e^{il\theta}, \qquad (1.9)$$

where $x/|x| = (\cos(\theta), \sin(\theta))$.

The Rayleigh Conjecture (RC) states that the series (1.9) converges up to the boundary $\Gamma$. This conjecture is false for many domains, although it holds for some domains, e.g., for a ball, see [3], [17], [23]. Recently A.G. Ramm [27] established a Modified Rayleigh Conjecture (MRC). A numerical implementation of the MRC method for obstacle scattering problems is presented in [28].

## 2. Modified Rayleigh Conjecture (MRC)

A 2-D version of the main result from [27], [28] is

**Theorem 2.1** *Let (1.9) be the unique representation of the outgoing solution $v(x)$ of the Exterior Dirichlet Problem (1.1)-(1.3). Fix*

an $\varepsilon > 0$.

Then there exists a positive integer $L = L(\varepsilon)$ and the coefficients $c_l = c_l(\varepsilon)$, $l = 0, \pm 1, ..., \pm L$ such that

(i). $\left\| f(x) - \sum_{l=-L}^{L} c_l H_l^{(1)}(k|x|)e^{il\theta} \right\|_{L^2(\Gamma)} \leq \varepsilon$,

(ii). $\left\| v(x) - \sum_{l=-L}^{L} c_l H_l^{(1)}(k|x|)e^{il\theta} \right\| = O(\varepsilon)$, $\varepsilon \to 0$,

where

$$\|\cdot\| = \|\cdot\|_{H_{loc}^m(D')} + \|\cdot\|_{L^2(D';(1+|x|)^{-\gamma})},$$

$\gamma > 1$, $m > 0$ is an arbitrary integer, $H^m$ is the Sobolev space, and

(iii). $c_l(\varepsilon) \to a_l$, as $\varepsilon \to 0$, $l = 0, \pm 1, \pm 2, ....$

According to Theorem 2.1 the computation of the outgoing solution of (1.1)-(1.3) is reduced to the approximation of the boundary values in (1.3). A direct implementation of the above algorithm is efficient for domains $D$ not very different from a circle, e.g. for an ellipse with a small eccentricity, but it fails for more complicated regions. The numerical difficulties happen because the Neumann functions $N_l$ with large values of $l$ are bigger than $N_l$ with small values of $l$ by many orders of magnitude. A finite precision of numerical computations makes it necessary to keep the values of $L$ not too high, e.g. $L \leq 20$. This restriction can be remedied by the following modification of the above algorithm, see [28]:

**Theorem 2.2** *Let $v(x)$ be the outgoing solution of the Exterior Dirichlet Problem (1.1)-(1.3). Suppose that points $x_1, x_2, ..., x_J$ are in the interior of $D$, and $\varepsilon > 0$.*
*Then*
*(i). There exists a positive integer $L = L(\varepsilon)$ and the coefficients*
$c_{lj} = c_{lj}(\varepsilon)$, $l = 0, \pm 1, ..., \pm L$,
$j = 1, 2, ..., J$ *such that*

$$\left\| f(x) - \sum_{j=1}^{J}\sum_{l=-L}^{L} c_{lj} H_l^{(1)}(k|x-x_j|)e^{il\theta_j} \right\|_{L^2(\Gamma)} \leq \varepsilon,$$
(2.1)

*where* $(x - x_j)/|x - x_j| = e^{i\theta_j}$.

*(ii) Let*
$$v_\varepsilon(x) = \sum_{j=1}^{J}\sum_{l=-L}^{L} c_{lj} H_l^{(1)}(k|x-x_j|)e^{il\theta_j}$$
(2.2)
*then*
$$\|v(x) - v_\varepsilon(x)\| = O(\varepsilon), \quad \varepsilon \to 0,$$

*where* $\|\cdot\| = \|\cdot\|_{H_{loc}^m(D')} + \|\cdot\|_{L^2(D';(1+|x|)^{-\gamma})}$,
$\gamma > 1$, $m > 0$ *is an arbitrary integer, $H^m$ is the Sobolev space.*

*(iii) The far field pattern of the approximate solution $v_\varepsilon(x)$ is given by*

$$A_{v_\varepsilon}(\alpha') = \sqrt{\frac{2}{\pi k}} e^{-i\frac{\pi}{4}} \sum_{j=1}^{J} e^{-ik\alpha' \cdot x_j} \sum_{i=-L}^{L} c_{lj}(-i)^l e^{il\theta}$$
(2.3)

*where* $\alpha' = x/|x| = e^{i\theta}$.

## 3. Direct Scattering Problem via MRC

According to Theorem 2.2 one can approximate the scattered field $u^s$ of the Direct Obstacle Scattering problem (1.4)-(1.6) by minimizing (2.1) with $f(x) = -u^i(x)$, $x \in \Gamma$. More precisely, the algorithm can be described as follows:

**Initialization.** Fix an integer $L > 0$ and an $\varepsilon > 0$. Choose $x_1, x_2, ..., x_J$ in the interior of $D$. If $\mathbf{r}(t)$, $0 \leq t < 2\pi$ is an equation of the boundary $\Gamma$, let

$$\psi_{lj}(t) = H_l^{(1)}(k\,|\mathbf{r}(t) - x_j|)e^{il\theta_j(t)},$$
$$j = 1,\,2,...,\,J,\quad l = 0, \pm 1, \pm 2,..., \pm L\,, \quad (3.1)$$
$$g(t) = -u^i(\mathbf{r}(t))) = -e^{ik\mathbf{r}(t)\cdot\alpha}\,, \quad (3.2)$$

where $(\mathbf{r}(t) - x_j)/|\mathbf{r}(t) - x_j| = e^{i\theta_j(t)}$.

**Minimization.** Minimize

$$\Phi(\mathbf{c}) = \left\| g(t) - \sum_{j=1}^{J}\sum_{l=-L}^{L} c_{lj}\psi_{lj}(t) \right\|_{L^2(0,\,2\pi)},\quad \mathbf{c} = \{c_{lj}\}\,,$$
$$(3.3)$$

for $\mathbf{c} \in \mathbb{C}^N$, $N = (2L+1)J$.

If the minimum of $\Phi$ in (3.3) is smaller than the prescribed tolerance $\varepsilon$, then the scattered field is approximated by $v_\varepsilon(x), x \in D'$, given by (2.2), and its Scattering Amplitude $A(\alpha',\,\alpha)$ is computed by formula (2.3).

The numerical implementation of the minimization algorithm begins with the choice of $M$ knots $0 = t_1 < t_2 < \cdots < t_M < 2\pi$, and points $x_j,\,j = 1,...,\,J$ in the interior of $D$. Then the values $\{\psi_{lj}(t_m)\}_{m=1}^{M}$ form $N = (2L+1)J$ vectors $\mathbf{a}^{(n)}$, $n = 1,\,2,...,\,N$ of length $M$. Let $\mathbf{b} = \{-u^i(t_m)\}_{m=1}^{M}$. Then the minimization problem (3.3) is reduced to the finite dimensional minimization problem

$$\min\{\|A\mathbf{c} - \mathbf{b}\|,\ \mathbf{c} \in \mathbb{C}^N\}\,,\quad (3.4)$$

where $A$ is the matrix containing vectors $\mathbf{a}^{(n)}$, $n = 1,\,2,...,\,N$ as its columns. If other outgoing solutions are used in addition to functions $\psi_{lj}$, the size of matrix $A$ is increased accordingly.

We use the Singular Value Decomposition (SVD) method (see e.g. [22]) to minimize (3.4). Small singular values of the matrix $A$ are used to identify and delete linearly dependent or almost linearly dependent combinations of vectors $\mathbf{a}^{(n)}$. This spectral cut-off makes the minimization process stable. The entire algorithm is summarized below. We denote by $V^H$ the complex conjugate transpose of a matrix $V$. Also, by the definition, the inner product in $\mathbb{C}^N$ complex conjugates its first component.

**Iterative MRC.** Fix an $\varepsilon > 0$, an integer $L > 0$, and $w_{\min} > 0$. Choose $M$ knots $0 = t_1 < t_2 < \cdots < t_M < 2\pi$, and points $x_j,\,j = 1,...,\,J$ in the interior of $D$.
Let $N = (2L+1)J$.

**(1) Initialization.**

a) Form vectors

$$\mathbf{a}^{(n)} = \{\psi_{lj}(t_m)\}_{m=1}^{M}\,,\ |l| \leq L\,,\ j = 1,\,2,...,\,J\,,$$

and the matrix $A$ of size $M \times N$, whose columns are the vectors $\mathbf{a}^{(n)}$.

b) Form vector

$$\mathbf{b} = \{-u^i(t_m)\}_{m=1}^{M}\,.$$

c) Use Singular Value Decomposition method to represent matrix $A$ as

$$A = UWV^H\,,$$

where the $M \times N$ matrix $U$ has orthonormal columns $\mathbf{u}^{(n)}$, $n = 1,...,\,N$, the square $N \times N$ matrix $V$ has orthonormal columns $\mathbf{v}^{(n)}$, $n = 1,...,\,N$, and the diagonal square $N \times N$ matrix $W = (w_n)_{n=1}^{N}$ is composed of the (nonnegative) singular values of $A$.

d) Let $\Sigma \subset \{1,\,2,...,\,N\}$ be defined by

$$\Sigma = \{n : w_n \geq w_{\min}\}\,.$$

e) Let $p = 0$.

(2) **Iterative step.**

a) Let $p := p+1$.
b) Form the set
   $P = \{n \in \Sigma : w_n \text{ is among } p \text{ largest singular values of } A\}$.
c) Compute the normalized residual
   $$r_p^{\min} = \frac{1}{\sqrt{M}} \sqrt{\|\mathbf{b}\|^2 - \sum_{n \in P} |<\mathbf{u}^{(n)}, \mathbf{b}>|^2}.$$

**(3) Stopping criterion.**
a) If $r_p^{\min} \leq \varepsilon$, then stop.
   The minimizer is given by
   $$\mathbf{c} = \sum_{n \in P} \frac{1}{w_n} <\mathbf{u}^{(n)}, \mathbf{b}> \mathbf{v}^{(n)}.$$
   Compute the scattered field $v_\varepsilon^s$ using (2.2) and the Far Field Pattern using (2.3).
b) If $r_p^{\min} > \varepsilon$, and $P \neq \Sigma$ repeat the iterative step (2).
c) Otherwise, the procedure has failed.

## 4. Numerical Experiments for MRC for obstacle scattering problems

TABLE 1. Normalized residuals attained in the numerical experiments.

| Experiment | $J$ | $k$ | $\alpha$ | $r^{\min}$ |
|---|---|---|---|---|
| I   | 4  | 1.0 | (1.0,0.0) | 0.000201 |
|     | 4  | 1.0 | (0.0,1.0) | 0.000357 |
|     | 4  | 5.0 | (0.0,0.0) | 0.001309 |
|     | 4  | 5.0 | (0.0,1.0) | 0.007228 |
| II  | 16 | 1.0 | (1.0,0.0) | 0.003555 |
|     | 16 | 1.0 | (0.0,1.0) | 0.002169 |
|     | 16 | 5.0 | (0.0,0.0) | 0.009673 |
|     | 16 | 5.0 | (0.0,1.0) | 0.007291 |
| III | 16 | 1.0 | (1.0,0.0) | 0.008281 |
|     | 16 | 1.0 | (0.0,1.0) | 0.007523 |
|     | 16 | 5.0 | (0.0,0.0) | 0.021571 |
|     | 16 | 5.0 | (0.0,1.0) | 0.024360 |
| IV  | 32 | 1.0 | (1.0,0.0) | 0.006610 |
|     | 32 | 1.0 | (0.0,1.0) | 0.006785 |
|     | 32 | 5.0 | (0.0,0.0) | 0.034027 |
|     | 32 | 5.0 | (0.0,1.0) | 0.040129 |

The results obtained by the MRC method (for smooth boundary $\Gamma$) were compared to the results obtained by the Boundary Integral Equation Method (BIEM) presented in [10]. Accordingly, to find the outgoing solution $v(x)$ of the Exterior Dirichlet Problem (1.1)-(1.3) one has to solve the integral equation

$$\varphi(x) + \int_\Gamma \left\{ \frac{\partial \Phi(x,y)}{\partial \nu(y)} - i\eta \Phi(x,y) \right\} \varphi(y) ds(y) = 2f(x), \quad x \in D' \quad (4.1)$$

for the density function $\varphi \in C(\Gamma)$. In the above equation $\Phi := \frac{i}{4} H_0^1(k|x-y|)$, and $H_0^1(r)$ is the Hankel function. Following the recommendations in [10] the value for the real coupling parameter $\eta$ was chosen to be equal to the wave number $k$. The above integral equation was solved using the Nyström method, see [10], Section 3.5.

After the density $\varphi$ is computed, the Far Field Pattern can be obtained from

$$A_\nu(\alpha') = \frac{e^{-i\frac{\pi}{4}}}{\sqrt{8\pi k}} \int_\Gamma \{k\nu(y) \cdot \alpha' + \eta\} e^{-ik\alpha' \cdot y} \varphi(y) ds(y)$$

(4.2)

where $\nu$ is the exterior unit normal vector to the boundary $\Gamma$.

We conducted numerical experiments for four obstacles: two ellipses of different eccentricity, a kite, and a triangle. Each case was tested for wave numbers $k = 1.0$ and $k = 5.0$. Each obstacle was subjected to incident waves corresponding to $\alpha = (1.0, 0.0)$ and $\alpha = (0.0, 1.0)$. The results are shown in Table 1. The column $J$ shows the number of the interior points $x_j$ used in the approximation (3.3). The choice of the points $x_j$ was different in each experiment. It is indicated below together with the description of the experiments. The column $r^{\min}$ shows the smallest value of the normalized residual achieved by the MRC minimization method in step 2.c of the Iterative MRC algorithm described in the previous section. Values $L = 5$ and $M = 720$ were used in all the experiments. The 720 knots $t_m$ were uniformly distributed on $[0, 2\pi]$. Under these conditions the relative error between the MRC and BIEM (implemented as described in the beginning of this section with $n = 64$ (see [10]) approximations of the scattering amplitude was less then 0.01%. See [28] for details.

FIGURE 1. Obstacle used in experiment II.

**Experiment I.** The boundary $\Gamma$ is the ellipse described by

$$\mathbf{r}(t) = (2.0\cos t, \ \sin t), \ 0 \le t < 2\pi \ . \quad (4.3)$$

The MRC minimization used 4 interior points $x_j = 0.7\mathbf{r}(\frac{\pi(j-1)}{2}), \ j = 1,..., 4$. Run time for the MRC was 2 seconds vs. 25 seconds for the BIEM on a 333 MHz PC.

**Experiment II.** The kite-shaped boundary $\Gamma$ (see [10], Section 3.5) is described by

$$\mathbf{r}(t) = (-0.65 + \cos t + 0.65\cos 2t, \ 1.5\sin t),$$
$$0 \le t < 2\pi \ , \quad (4.4)$$

See Figure 1. The MRC minimization used 16 interior points $x_j = 0.9\mathbf{r}(\frac{\pi(j-1)}{8}), \ j = 1,..., 16$. Run time for the MRC was 33 seconds vs. 44 seconds for the BIEM.

**Experiment III.** The boundary $\Gamma$ is the triangle with vertices at $(-1.0, \ 0.0)$ and $(1.0, \pm 1.0)$. The MRC minimization used 16 interior points $x_j = 0.9\mathbf{r}(\frac{\pi(j-1)}{8}), \ j = 1,..., 16$. Run time for the MRC was about 30 seconds.

**Experiment IV.** The boundary $\Gamma$ is the ellipse described by

$$\mathbf{r}(t) = (0.1\cos t, \ \sin t), \ 0 \le t < 2\pi \ . \quad (4.5)$$

The MRC minimization used 32 interior points $x_j = 0.95\mathbf{r}(\frac{\pi(j-1)}{16}), \ j = 1,..., 32$. Run time for the MRC was about 140 seconds.

## 5. Scattering by periodic structures

Determination of fields scattered by periodic structures is of a great importance in modern diffractive optics, and there is a vast literature on both the direct and inverse problems of this type, see, for example [21]. Still, an efficient computation of such fields presents certain difficulties. In the next Sections we present some theoretical background, a modification of the MRC method, and numerical results for such a scattering, see [29].

For simplicity we consider a 2-D setting, but our arguments can be as easily applied to $n$-dimensional problems, $n \ge 2$. Let $f : \mathbb{R} \to \mathbb{R}$, $f(x+L) = f(x)$ be an $L$-periodic Lipschitz continuous function, and let $D$ be the domain

$$D = \{(x, y): y \ge f(x), x \in \mathbb{R}\}$$

Without loss of generality we assume that $f \ge 0$. If it is not, one can choose the origin so that this assumption is satisfied, because $M := \sup_{0 \le x \le L} |f(x)| < \infty$.

Let $\mathbf{x} = (x, \ y)$ and $u(\mathbf{x})$ be the total field satisfying

$$(\Delta + k^2)u = 0, \ \mathbf{x} \in D, \ k = const > 0 \quad (5.1)$$

$$u = 0 \text{ on } S := \partial D, \quad (5.2)$$

$$u = u_0 + v, \ u_0 := e^{ik\alpha \cdot \mathbf{x}}, \quad (5.3)$$

where the unit vector $\alpha = (\cos\theta, -\sin\theta)$, $0 < \theta < \pi/2$, and $v(\mathbf{x})$ is the scattered field, whose asymptotic behavior as $y \to \infty$ will be specified below, and

$$u(x+L, y) = \nu u(x, y),$$
$$u_x(x+L, y) = \nu u_x(x, y) \quad (5.4)$$
$$\text{in } D, \quad \nu := e^{ikL\cos\theta}.$$

Conditions (5.4) are the $qp$ (**quasiperiodicity**) conditions. To find the proper radiation condition for the scattered field $v(\mathbf{x})$ consider the spectral problem

$$\varphi'' + \lambda^2 \varphi = 0, \quad 0 < x < L, \quad (5.5)$$

$$\varphi(L) = \nu\varphi(0), \quad \varphi'(L) = \nu\varphi'(0) \quad (5.6)$$

arising from the separation of variables in (5.1)-(5.4). This problem has a discrete spectrum, and its eigenfunctions form a basis in $L^2(0, L)$. One has

$$\varphi = Ae^{i\lambda x} + Be^{-i\lambda x}, \quad A, B = const,$$

$$Ae^{i\lambda L} + Be^{-i\lambda L} = \nu(A+B),$$
$$i\lambda Ae^{i\lambda L} - i\lambda Be^{-i\lambda L} = i\lambda\nu(A-B)$$

Thus

$$\begin{vmatrix} e^{i\lambda L} - \nu & e^{-i\lambda L} - \nu \\ i\lambda(e^{i\lambda L} - \nu) & -i\lambda(e^{-i\lambda L} - \nu) \end{vmatrix} = 0.$$

So

$$i\lambda(e^{i\lambda L} - \nu)(e^{-i\lambda L} - \nu) = 0.$$

If $\lambda = 0$, then $\varphi = A + Bx$, $A + BL = \nu A$, $B = \nu B$. Since $\nu = e^{ikL\cos\theta}$, one has no eigenvalue $\lambda = 0$ unless $kL\sin\theta = 2\pi m$, $m > 0$ is an integer. Let us assume that $kL\cos\theta \neq 2\pi m$. Then

$$e^{i\lambda L} = e^{ikL\cos\theta} \quad \text{or} \quad e^{-i\lambda L} = e^{ikL\cos\theta},$$

that is

$$\lambda_j^+ = k\cos\theta + \frac{2\pi j}{L}, \quad \text{or} \quad \lambda_j^- = -k\cos\theta + \frac{2\pi j}{L},$$
$$j = 0, \pm 1, \pm 2,...$$

The corresponding eigenfunctions are $e^{i\lambda_j^+ x}$ and $e^{-i\lambda_j^+ x}$. We will use the system $e^{i\lambda_j^+ x}$, which forms an orthogonal basis in $L^2(0, L)$. One has:

$$\int_0^L e^{i\lambda_j^+ x} e^{-i\lambda_m^+ x} \, dx = \int_0^L e^{\frac{2\pi i}{L}(j-m)} \, dx = 0, \quad j \neq m.$$

The normalized eigenfunctions are

$$\varphi_j(x) = \frac{e^{i\lambda_j^+ x}}{\sqrt{L}}, \quad j = 0, \pm 1, \pm 2,...$$

These functions form an orthonormal basis of $L^2(0, L)$. Let us look for $v(\mathbf{x}) = v(x, y)$ of the form

$$v(x, y) = \sum_{j=-\infty}^{\infty} c_j v_j(y)\varphi_j(x), \quad y > M,$$
$$c_j = const \quad (5.7)$$

For $y > M$, equation (5.1) implies

$$v_j'' + (k^2 - \lambda_j^2)v_j = 0. \quad (5.8)$$

Let us assume that $\lambda_j^2 \neq k^2$ for all $j$. Then

$$v_j(y) = e^{i\mu_j y}, \quad (5.9)$$

where, for finitely many $j$, the set of which is denoted by $J$, one has:

$$\mu_j = (k^2 - \lambda_j^2)^{1/2} > 0, \quad \text{if } \lambda_j^2 < k^2, \, j \in J,$$
$$(5.10)$$

and

$$\mu_j = i(\lambda_j^2 - k^2)^{1/2}, \text{ if } \lambda_j^2 > k^2, \quad j \notin J. \tag{5.11}$$

The **radiation condition** at infinity requires that the scattered field $v(x, y)$ be representable in the form (5.7) with $v_j(y)$ defined by (5.9)-(5.11).

The **Periodic Scattering Problem** consists of finding the solution to (5.1)-(5.4) satisfying the radiation condition (5.7), (5.9)-(5.11).

The existence and uniqueness for such a scattering problem is established in [29]. In [1] the scattering by a periodic structure was considered earlier, and was based on a uniqueness theorem from [11]. There are many papers on scattering by periodic structures, of which we mention a few [1], [2], [4], [5], [7], [12], [13], [18], [19], [21], [46]. The Rayleigh conjecture is discussed in several of the above papers. It was shown (e.g. [21], [3]) that this conjecture is incorrect, in general. As we have already discussed in the previous sections, the modified Rayleigh conjecture is a theorem proved in [27] for scattering by bounded obstacles.

Here we present an outline of one of the existence and uniqueness proofs for scattering by periodic structures following [29]. To prove the existence and uniqueness of the resolvent kernel $G$ we start by constructing a scattering theory quite similar to the one for the exterior of a bounded obstacle [23]. The *first step* is to construct an analog to the half-space Dirichlet Green's function. The function $g = g(\mathbf{x}, \boldsymbol{\xi}, k)$ can be constructed analytically ($\mathbf{x} = (x_1, x_2)$, $\boldsymbol{\xi} = (\xi_1, \xi_2)$):

$$g(\mathbf{x}, \boldsymbol{\xi}) = \sum_j \varphi_j(x_1)\overline{\varphi_j(\xi_1)}g_j(x_2, \xi_2, k) \tag{5.12}$$

$$g_j := g_j(x_2, \xi_2, k) = \begin{cases} v_j(x_2)\psi_j(\xi_2), & x_2 > \xi_2 \\ v_j(\xi_2)\psi_j(x_2), & x_2 < \xi_2 \end{cases}$$

$$\psi_j = (\mu_j)^{-1} e^{i\mu_j b} \sin[\mu_j(\xi_2 + b)],$$

$$\mu_j = [k^2 - \lambda_j^2]^{1/2}, \quad v_j(x_2) = e^{i\mu_j x_2},$$

where

$$\psi_j'' + (k^2 - \lambda_j^2)\psi_j = 0, \quad \psi_j(-b) = 0,$$

$$W[v_j, \psi_j] = 1, \quad \lambda_j = k\cos(\theta) + \frac{2\pi j}{L},$$

and $W[v, \psi]$ is the Wronskian.

The function $g$ is analytic with respect to $k$ on the complex plain with cuts along the rays $\lambda_j - i\tau$, $0 \le \tau < \infty$, $j = 0, \pm 1, \pm 2, ...$, in particular, in the region $\Im k > 0$, up to the real positive half-axis except for the set $\{\lambda_j\}_{j=0,\pm 1,\pm 2,...}$.

Choose $b > 0$ such that $k^2 > 0$ is not an eigenvalue of the problem:

$$(\Delta + k^2)\psi = 0,$$
$$\text{in } D_{-b} := \{(x, y) : -b \le y \le f(x), \ 0 \le x \le L\} \tag{5.13}$$

$$\psi|_{y=-b} = 0, \ \psi_N = 0 \text{ on } S,$$
$$\psi(x+L, y) = \nu\psi(x, y),$$
$$\psi_x(x+L, y) = \nu\psi_x(x, y) \tag{5.14}$$

One has
$$(\Delta + k^2)g = -\delta(\mathbf{x}-\boldsymbol{\xi}), \ \mathbf{x} = (x_1, x_2),$$
$$\boldsymbol{\xi} = (\xi_1, \xi_2), \tag{5.15}$$
$$\mathbf{x} \in \{(x, y) : -b < y < \infty, \ 0 \le x \le L\},$$

$$g|_{y=-b} = 0, \tag{5.16}$$

and

$$(\Delta + k^2)G = -\delta(\mathbf{x}-\boldsymbol{\xi}), \tag{5.17}$$

$G$ satisfies the $qp$ condition and the radiation condition (it is outgoing at infinity).

Multiply (5.15) by $G$, (5.17) by $g$, subtract from the second equation the first one, integrate over $D_{LR}$, and take $R \to \infty$, to get

$$G = g + \int_{S_L} (Gg_N - G_N g) ds = g - \int_{S_L} g\mu\ ds,$$
$$\mu := G_N |_{S_L}.$$
(5.18)

The $qp$ condition allows one to cancel the integrals over the lateral boundary ($x=0$ and $x=L$), and the radiation condition allows one to have

$$\lim_{R\to\infty}(Gg_N - G_N g) ds = 0.$$

Differentiate (5.18) to get

$$\mu = -A\mu + 2\frac{\partial g}{\partial N} \text{ on } S_L,$$
$$A\mu := 2\int_{S_L} \frac{\partial g(s,\ \sigma)}{\partial N_s} \mu(\sigma)\ d\sigma.$$
(5.19)

This is a Fredholm equation for $\mu$ in $L^2(S_L)$, if $S_L$ is $C^{1,m}$, $m>0$. The homogeneous equation (5.19) has only the trivial solution: if $\mu + A\mu = 0$, then the function $\psi := \int_{S_L} g\mu\ ds$ satisfies $\psi_N^+|_{S_L} = 0$, where $\psi_N^+(\psi_N^-)$ is the normal derivative of $\psi$ from $D_{-b}(D_L)$, and we use the known formula for the normal derivative of the single layer potential at the boundary. The $\psi$ satisfies also (5.13) and (5.14), and, by the choice of $b$, one has $\psi = 0$ in $D_{-b}$. Also $\psi = 0$ in $D_L$, because $(\Delta + k^2)\psi = 0$ in $D_L$, $\psi|_{S_L} = 0$ (by the continuity of the single layer potential), $\psi$ satisfies the $qp$ condition (because $g$ satisfies it), and $\psi$ is outgoing (because $g$ is).

Since $\psi = 0$ in $D_{-b}$ and in $D_L$, one concludes that $\mu = \psi_N^+ - \psi_N^-$, where $\psi_N^+(\psi_N^-)$ is the normal derivative of $\psi$ from $D_{-b}(D_L)$, and we use the jump relation for the normal derivative of the single layer potential.

Thus, we have proved the existence and uniqueness of $\mu$, and, therefore, of $G$, and got a representation formula

$$G = g - \int_{S_L} g\mu\ ds. \qquad (5.20)$$

## 6. MRC for scattering by periodic structures

Rayleigh conjectured [46] ("Rayleigh hypothesis") that the series (5.7) converges up to the boundary $S_L$. This conjecture is wrong (see [21]) for some $f(x)$. Since the Rayleigh hypothesis has been widely used for numerical solution of the scattering problem by physicists and engineers, and because these practitioners reported high instability of the numerical solution, and there are no error estimates, we propose a modification of the Rayleigh conjecture, which is a Theorem. This MRC (Modified Rayleigh Conjecture) can be used for a numerical solution of the scattering problem, and it gives an error estimate for this solution. Our arguments are very similar to the ones in [27].

Rewrite the scattering problem (5.1)-(5.4) as

$$(\Delta + k^2)v = 0 \text{ in } D,\ v = -u_0 \text{ on } S_L, (6.1)$$

where $v$ satisfies (5.4), and $v$ has representation (5.7), that is, $v$ is "outgoing", it satisfies the radiation condition. Fix an arbitrarily small $\varepsilon > 0$, and assume that

$$\left\| u_0 + \sum_{|j| \leq j(\varepsilon)} c_j(\varepsilon) v_j(y) \varphi_j(x) \right\| \leq \varepsilon,$$
$$0 \leq x \leq L,\ y = f(x),$$
(6.2)

where $\|\cdot\| = \|\cdot\|_{L^2(S_L)}$.

**Lemma 6.1** *For any $\varepsilon > 0$, however small, and for any $u_0 \in L^2(S_L)$, there exists $j(\varepsilon)$ and $c_j(\varepsilon)$ such that (6.2) holds.*

*Proof.* Let us prove the completeness of the system $\{\varphi_j(x) v_j(f(x))\}_{j=0,\pm 1,\pm 2,...}$ in $L^2(S_L)$. Assume that there is an $h \in L^2(S_L)$, $h \not\equiv 0$

such that

$$\int_{S_L} h\overline{\varphi_j(x)} v_j(f(x))\, ds = 0 \qquad (6.3)$$

for any $j$. From (6.3) one derives (cf. [23], p.162-163)

$$\psi(\mathbf{x}) := \int_{S_L} hg(\mathbf{x}, \boldsymbol{\xi}) d\xi = 0, \quad \mathbf{x} \in D_{-b}. \quad (6.4)$$

Thus $\psi = 0$ in $D_L$, and $h = \psi_N^+ - \psi_N^- = 0$. Lemma is proved.

**Lemma 6.2.** *If (6.2) holds, then*

$$\left\| v(\mathbf{x}) - \sum_{|j| \le j(\varepsilon)} c_j(\varepsilon) v_j(y) \varphi_j(x) \right\| \le c\varepsilon,$$
$$\forall x,\ y \in D_L,\ 0 \le x \le L,\ y = f(x),$$

*where $c = const > 0$ does not depend on $\varepsilon,\ x,\ y$; $R > M$ is an arbitrary fixed number, and*

$$\|w\| = \sup_{\mathbf{x} \in D \setminus D_{LR}} |w(\mathbf{x})| + \|w\|_{H^{1/2}(D_{LR})}.$$

*Proof.* Let $w := v - \sum_{|j| \le j(\varepsilon)} c_j(\varepsilon) v_j(y) \varphi_j(x)$. Then $w$ solves equation (5.1), $w$ satisfies (5.4), $w$ is outgoing, and $\|w\|_{L^2(S_L)} \le \varepsilon$. One has

$$w(\mathbf{x}) = -\int_{S_L} w G_N(\mathbf{x}, \boldsymbol{\xi})\, ds. \qquad (6.5)$$

Thus (6.2), i.e. $\|w\| := \|w\|_{L^2(S_L)} \le \varepsilon$, implies

$$|w(\mathbf{x})|_{y=R} \le \|w\|_{L^2(S_L)} \|G_N(\mathbf{x}, \boldsymbol{\xi})\|_{L^2(S_L)} \le c\varepsilon,$$
$$c = const > 0,$$
$$(6.6)$$

where $c$ is independent of $\varepsilon$, and $R > \max f(x)$ is arbitrary. Now let us use the elliptic inequality

$$\|w\|_{H^m(D_{LR})} \le c\left( \|w\|_{H^{m-0.5}(S_L)} + \|w\|_{H^{m-0.5}(S_R)} \right),$$
$$(6.7)$$

where we have used the equation $\Delta w + k^2 w = 0$, and assumed that $k^2$ is not a Dirichlet eigenvalue of the Laplacian in $D_{LR}$, which can be done without loss of generality, because one can vary $R$. The integer $m \ge 0$ is arbitrary if $S_L$ is sufficiently smooth, and $m \le 1$ if $S_L$ is Lipschitz. Taking $m = 0.5$ and using (6.2) and (6.6) one gets

$$\|w\|_{H^{1/2}(D_{LR})} \le c\varepsilon. \qquad (6.8)$$

Thus, in a neighborhood of $S_L$, we have proved estimate (6.8), and in the complement of this neighborhood in $D_L$ we have proved estimate (6.6). Lemma is proved.

*Remark 6.3.* In (6.7) there are no terms with boundary norms over the lateral boundary (lines $x = 0$ and $x = L$) because of the quasiperiodicity condition.

From Lemma 6.2 the basic result, Theorem 6.4, follows immediately:

**Theorem 6.4. MRC-Modified Rayleigh Conjecture.**

*Fix $\varepsilon > 0$, however small, and choose a positive integer $p$. Find*

$$\min_{c_j} \left\| u_0 + \sum_{|j| \le p} c_j \varphi_j(x) v_j(y) \right\| := m(p). \qquad (6.9)$$

*Let $\{c_j(p)\}$ be the minimizer of (6.9). If $m(p) \le \varepsilon$, then*

$$v(p) = \sum_{|j| \le p} c_j(p) \varphi_j(x) v_j(y) \qquad (6.10)$$

*satisfies the inequality*

$$\|v - v(p)\| \le c\varepsilon, \qquad (6.11)$$

where $c = const > 0$ does not depend on $\varepsilon$. If $m(p) > \varepsilon$, then there exists $j = j(\varepsilon) > p$ such that $m(j(\varepsilon)) < \varepsilon$. Denote $c_j(j(\varepsilon)) := c_j(\varepsilon)$ and $v(j(\varepsilon)) := v_\varepsilon$. Then

$$\|v - v_\varepsilon\| \leq c\varepsilon. \quad (6.12)$$

## 7. Numerical solution of the periodic scattering problem

According to the MRC method (Theorem 6.4), if the restriction of the incident field $-u_0(x, y)$ to $S_L$ is approximated as in (6.9), then the series (6.10) approximates the scattered field in the entire region above the profile $y = f(x)$. However, a numerical method that uses (6.9) does not produce satisfactory results as reported in [21] and elsewhere. Our own numerical experiments confirm this observation. A way to overcome this difficulty is to realize that the numerical approximation of the field $-u_0 |_{S_L}$ can be carried out by using outgoing solutions described below.

Let $\boldsymbol{\xi} = (\xi_1, \xi_2) \in D_{-b}$, where $b > 0$,
$$D_{-b} := \{(\xi_1, \xi_2) : -b \leq \xi_2 \leq f(x), \ 0 \leq \xi_1 \leq L\},$$
and $g(\mathbf{x}, \boldsymbol{\xi})$ be defined as in Section 5. Then $g(\mathbf{x}, \boldsymbol{\xi})$ is an outgoing solution satisfying $\Delta g + k^2 g = 0$ in $D_L$, according to (5.15).

To implement the MRC method numerically one proceeds as follows:

1) Choose the nodes $\mathbf{x}_i, \ i = 1, 2, ..., N$ on the profile $S_L$. These points are used to compute approximate $L^2$ norms on $S_L$.
2) Choose points $\boldsymbol{\xi}^{(1)}, \boldsymbol{\xi}^{(2)}, ..., \boldsymbol{\xi}^{(M)}$ in $D_{-b}$, $M < N$.
3) Form the vectors $\mathbf{b} = (u_0(\mathbf{x}_i))$, and $\mathbf{a}^{(m)} = (g(\mathbf{x}_i, \boldsymbol{\xi}^{(m)}))$, $i = 1, 2, ..., N$, $m = 1, 2, ..., M$. Let $\mathbf{A}$ be the $N \times M$ matrix containing vector $\mathbf{a}^{(m)}$ as its columns.
4) Find the Singular Value Decomposition of $\mathbf{A}$. Use a predetermined $w_{\min} > 0$ to eliminate its small singular values. Use the decomposition to compute
$$r^{\min} = \min\{\|\mathbf{b} + \mathbf{Ac}\|, \ \mathbf{c} \in \mathbb{C}^M\}, \text{ where}$$
$$\|\mathbf{a}\|^2 = \frac{1}{N}\sum_{i=1}^{N} |a_i|^2.$$
5) **Stopping criterion.** Let $\varepsilon > 0$.
   a) If $r^{\min} \leq \varepsilon$, then stop. Use the coefficients $\mathbf{c} = \{c_1, c_2, ..., c_M\}$ obtained in the above minimization step to compute the scattered field by
   $$v(x, y) = \sum_{m=1}^{M} c_m g(x, y, \boldsymbol{\xi}^{(m)}).$$
   b) If $r^{\min} > \varepsilon$, then increase $N$, $M$ by the order of 2, readjust the location of points $\boldsymbol{\xi}^{(m)} \in D_{-b}$ as needed, and repeat the procedure.

We have conducted numerical experiments for four different profiles. In each case we used $L = \pi$, $k = 1.0$ and three values for the angle $\theta$. Table 2 shows the resulting residuals $r^{\min}$. Note that $\|\mathbf{b}\| = 1$. Thus, in all the considered cases, the MRC method achieved 0.04% to 2% accuracy of the approximation. Other parameters used in the experiments were chosen as follows: $N = 256$, $M = 64$, $w_{\min} = 10^{-8}$, $b = 1.2$. The value of $b > 0$, used in the definition of $g$, was chosen experimentally, but the dependency of $r^{\min}$ on $b$ was slight. The Singular Value Decomposition (SVD) is used in Step 4 since the vectors $\mathbf{a}^{(m)}$, $m = 1, 2, ..., M$ may be nearly linearly dependent, which leads to an instability in the determination of the minimizer $\mathbf{c}$. According to the SVD method this instability is eliminated by cutting off small singular values of the matrix $\mathbf{A}$, see e.g. [22] for details. The cut-off value $w_{\min} > 0$ was chosen experimentally. We used the truncated series (5.12) with $|j| \leq 120$ to compute functions $g(x, y, \boldsymbol{\xi})$. A typical run time on a 333 MHz PC was about $40s$ for each experiment.

The following is a description of the profiles

$y = f(x)$, the nodes $\mathbf{x}_i \in S_L$, and the poles $\xi^{(m)} \in D_{-b}$ used in the computation of $g(\mathbf{x}_i, \xi^{(m)})$ in Step 3. For example, in profile I the $x$-coordinates of the $N$ nodes $\mathbf{x}_i \in S_L$ are uniformly distributed on the interval $0 \le x \le L$. The poles $\xi^{(m)} \in D_{-b}$ were chosen as follows: every fourth node $\mathbf{x}_i$ was moved by a fixed amount $-0.1$ parallel to the $y$ axis, so it would be within the region $D_{-b}$. The location of the poles was chosen experimentally to give the smallest value of the residual $r^{\min}$.

TABLE 2. Residuals attained in the numerical experiments.

| Profile | $\theta$ | $r^{\min}$ |
| --- | --- | --- |
| I | $\pi/4$ | 0.000424 |
|   | $\pi/3$ | 0.000407 |
|   | $\pi/2$ | 0.000371 |
| II | $\pi/4$ | 0.001491 |
|    | $\pi/3$ | 0.001815 |
|    | $\pi/2$ | 0.002089 |
| III | $\pi/4$ | 0.009623 |
|     | $\pi/3$ | 0.011903 |
|     | $\pi/2$ | 0.013828 |
| IV | $\pi/4$ | 0.014398 |
|    | $\pi/3$ | 0.017648 |
|    | $\pi/2$ | 0.020451 |

**Profile I.** $f(x) = \sin(2x)$ for $0 \le x \le L$,
$t_i = iL/N$, $\mathbf{x}_i = (t_i, f(t_i))$, $i = 1, 2, ..., N$,
$\xi^{(m)} = (x_{4m}, y_{4m} - 0.1)$, $m = 1, 2, ..., M$.

**Profile II.** $f(x) = \sin(0.2x)$ for $0 \le x \le L$,
$t_i = iL/N$, $\mathbf{x}_i = (t_i, f(t_i))$, $i = 1, 2, ..., N$,
$\xi^{(m)} = (x_{4m}, y_{4m} - 0.1)$, $m = 1, 2, ..., M$.

**Profile III.** $f(x) = x$ for $0 \le x \le L/2$, $f(x) = L - x$ for $L/2 \le x \le L$, $t_i = iL/N$,
$\mathbf{x}_i = (t_i, f(t_i))$, $i = 1, 2, ..., N$,
$\xi^{(m)} = (x_{4m}, y_{4m} - 0.1)$, $m = 1, 2, ..., M$.

**Profile IV.** $f(x) = x$ for $0 \le x \le L$,
$t_i = 2iL/N$, $\mathbf{x}_i = (t_i, f(t_i))$, $i = 1, ..., N/2$,
$\mathbf{x}_i = (L, f(2(i - N/2)L/N)), i = N/2 + 1, ..., N$,
$\xi^{(m)} = (x_{4m} - 0.03, y_{4m} - 0.05)$, $m = 1, 2, ..., M$

In this profile $N/2$ nodes $\mathbf{x}_i$ are uniformly distributed on its slant part, and $N/2$ nodes are uniformly distributed on its vertical portion $x = L$.

The experiments show that the MRC method provides a competitive alternative to other methods for the computation of fields scattered from periodic structures. It is fast and inexpensive. The results depend on the number of the internal points $\xi^{(m)}$ and on their location. A similar MRC method for the computation of fields scattered by a bounded obstacle was presented in [28].

## 8. The Support Function Method (SFM)

The Inverse Scattering Problem consists of finding the obstacle $D$ from the Scattering Amplitude, or similarly observed data. The Support Function Method (SFM) originally developed in a 3-D setting in [23], pp 94-99, is used to approximately locate the obstacle $D$. The method is derived using a high-frequency approximation to the scattered field for smooth, strictly convex obstacles. It turns out that this inexpensive method also provides a good localization of obstacles in the resonance region of frequencies. One can restate the SFM in a 2-D setting as follows (see [42]).

Let $D \subset \mathbb{R}^2$ be a smooth and strictly convex obstacle with the boundary $\Gamma$. Let $\nu(\mathbf{y})$ be the unique outward unit normal vector to $\Gamma$ at $\mathbf{y} \in \Gamma$. Fix an incident direction $\alpha \in S^1$. Then the boundary $\Gamma$ can be decomposed into the following two parts:

$$\Gamma_+ = \{\mathbf{y} \in \Gamma : \nu(\mathbf{y}) \cdot \alpha < 0\},$$
$$\text{and } \Gamma_- = \{\mathbf{y} \in \Gamma : \nu(\mathbf{y}) \cdot \alpha \geq 0\}, \quad (8.1)$$

which are, correspondingly, the illuminated and the shadowed parts of the boundary for the chosen incident direction $\alpha$.

Give $\alpha \in S^1$, its **specular point** $\mathbf{s}_0(\alpha) \in \Gamma_+$ is defined from the condition:

$$\mathbf{s}_0(\alpha) \cdot \alpha = \min_{\mathbf{s} \in \Gamma_+} \mathbf{s} \cdot \alpha \quad (8.2)$$

Note that the equation of the tangent line to $\Gamma_+$ at $\mathbf{s}_0$ is

$$<x_1, x_2> \cdot \alpha = \mathbf{s}_0(\alpha) \cdot \alpha, \quad (8.3)$$

and

$$\nu(\mathbf{s}_0(\alpha)) = -\alpha. \quad (8.4)$$

The **Support function** $d(\alpha)$ is defined by

$$d(\alpha) = \mathbf{s}_0(\alpha) \cdot \alpha. \quad (8.5)$$

Thus $|d(\alpha)|$ is the distance from the origin to the unique tangent line to $\Gamma_+$ perpendicular to the incident vector $\alpha$. Since the obstacle $D$ is assumed to be convex

$$D = \cap_{\alpha \in S^1} \{x \in \mathbb{R}^2 : x \cdot \alpha \geq d(\alpha)\} \quad (8.6)$$

The boundary $\Gamma$ of $D$ is smooth, hence so is the Support Function. The knowledge of this function allows one to reconstruct the boundary $\Gamma$ using the following procedure.

Parametrize unit vectors $\mathbf{l} \in S^1$ by $\mathbf{l}(t) = (\cos t, \sin t), \quad 0 \leq t < 2\pi$ and define

$$p(t) = d(\mathbf{l}(t)), \quad 0 \leq t < 2\pi. \quad (8.7)$$

Equation (8.3) and the definition of the Support Function give

$$x_1 \cos t + x_2 \sin t = p(t). \quad (8.8)$$

Since $\Gamma$ is the envelope of its tangent lines, its equation can be found from (8.8) and

$$-x_1 \sin t + x_2 \cos t = p'(t). \quad (8.9)$$

Therefore the parametric equations of the boundary $\Gamma$ are

$$x_1(t) = p(t)\cos t - p'(t)\sin t,$$
$$x_2(t) = p(t)\sin t + p'(t)\cos t. \quad (8.10)$$

So, the question is how to construct the Support function $d(\mathbf{l})$, $\mathbf{l} \in S^1$ from the knowledge of the Scattering Amplitude. In 2-D the Scattering Amplitude is related to the total field $u = u_0 + v$ by

$$A(\alpha', \alpha) = -\frac{e^{i\frac{\pi}{4}}}{\sqrt{8\pi k}} \int_\Gamma \frac{\partial u}{\partial \nu(\mathbf{y})} e^{-ik\alpha' \cdot \mathbf{y}} ds(\mathbf{y}). \quad (8.11)$$

In the case of the "soft" (i.e. the pressure field satisfies the Dirichlet boundary condition $u = 0$) the Kirchhoff (high frequency) approximation gives

$$\frac{\partial u}{\partial \nu} = 2\frac{\partial u_0}{\partial \nu} \quad (8.12)$$

on the illuminated part $\Gamma_+$ of the boundary $\Gamma$, and

$$\frac{\partial u}{\partial \nu} = 0 \quad (8.13)$$

on the shadowed part $\Gamma_-$. Therefore, in this approximation,

$$A(\alpha', \alpha) = -\frac{ike^{i\frac{\pi}{4}}}{\sqrt{2\pi k}} \int_{\Gamma_+} \alpha \cdot \nu(\mathbf{y}) e^{ik(\alpha-\alpha') \cdot \mathbf{y}} ds(\mathbf{y}). \quad (8.14)$$

Let $L$ be the length of $\Gamma_+$, and $\mathbf{y} = \mathbf{y}(\zeta)$, $0 \leq \zeta \leq L$ be its arc length parametrization. Then

$$A(\alpha', \alpha) = -\frac{i\sqrt{k} \, e^{i\frac{\pi}{4}}}{\sqrt{2\pi}} \int_0^L \alpha \cdot \nu(\mathbf{y}(\zeta)) e^{ik(\alpha-\alpha') \cdot \mathbf{y}(\zeta)} d\zeta. \quad (8.15)$$

Let $\zeta_0 \in [0, L]$ be such that $\mathbf{s}_0 = \mathbf{y}(\zeta_0)$ is the specular point of the unit vector $\mathbf{l}$, where

$$\mathbf{l} = \frac{\alpha - \alpha'}{|\alpha - \alpha'|} \quad . \tag{8.16}$$

Then $\nu(\mathbf{s}_0) = -\mathbf{l}$, and $d(\mathbf{l}) = \mathbf{y}(\zeta_0) \cdot \mathbf{l}$. Let

$$\varphi(\zeta) = (\alpha - \alpha') \cdot \mathbf{y}(\zeta) \quad .$$

Then $\varphi(\zeta) = \mathbf{l} \cdot \mathbf{y}(\zeta) |\alpha - \alpha'|$. Since $\nu(\mathbf{s}_0)$ and $\mathbf{y}'(\zeta_0)$ are orthogonal, one has

$$\varphi'(\zeta_0) = \mathbf{l} \cdot \mathbf{y}'(\zeta_0) |\alpha - \alpha'| = 0 \quad .$$

Therefore, due to the strict convexity of $D$, $\zeta_0$ is also the unique non-degenerate stationary point of $\varphi(\zeta)$ on the interval $[0, L]$, that is $\varphi'(\zeta_0) = 0$, and $\varphi''(\zeta_0) \neq 0$.

According to the Stationary Phase method

$$\int_0^L f(\zeta) e^{ik\varphi(\zeta)} d\zeta$$
$$= f(\zeta_0) \exp\left[ ik\varphi(\zeta_0) + \frac{i\pi}{4} \frac{\varphi''(\zeta_0)}{|\varphi''(\zeta_0)|} \right]$$
$$\sqrt{\frac{2\pi}{k|\varphi''(\zeta_0)|}} \left[ 1 + O(\frac{1}{k}) \right] \tag{8.17}$$

as $k \to \infty$.

By the definition of the curvature $\kappa(\zeta_0) = |\mathbf{y}''(\zeta_0)|$. Therefore, from the collinearity of $\mathbf{y}''(\zeta_0)$ and $\mathbf{l}$, $|\varphi''(\zeta_0)| = |\alpha - \alpha'| \kappa(\zeta_0)$. Finally, the strict convexity of $D$, and the definition of $\varphi(\zeta)$, imply that $\zeta_0$ is the unique point of minimum of $\varphi$ on $[0, L]$, and

$$\frac{\varphi''(\zeta_0)}{|\varphi''(\zeta_0)|} = 1 \quad . \tag{8.18}$$

Using (8.17)-(8.18), expression (8.15) becomes:

$$A(\alpha', \alpha) = -\frac{\mathbf{l} \cdot \alpha}{\sqrt{|\alpha - \alpha'| \kappa(\zeta_0)}} e^{ik(\alpha - \alpha') \cdot \mathbf{y}(\zeta_0)}$$
$$\left[ 1 + O\left(\frac{1}{k}\right) \right], \quad k \to \infty \quad . \tag{8.19}$$

At the specular point one has $\mathbf{l} \cdot \alpha' = -\mathbf{l} \cdot \alpha$. By the definition $\alpha - \alpha' = \mathbf{l} |\alpha - \alpha'|$. Hence $\mathbf{l} \cdot (\alpha - \alpha') = |\alpha - \alpha'|$ and $2\mathbf{l} \cdot \alpha = |\alpha - \alpha'|$. These equalities and $d(\mathbf{l}) = \mathbf{y}(\zeta_0) \cdot \mathbf{l}$ give

$$A(\alpha', \alpha) = \frac{1}{2} \sqrt{\frac{|\alpha - \alpha'|}{\kappa(\zeta_0)}} \; e^{ik|\alpha - \alpha'|d(\mathbf{l})}$$
$$\left[ 1 + O(\frac{1}{k}) \right], k \to \infty \tag{8.20}$$

Thus, the approximation

$$A(\alpha', \alpha) \approx -\frac{1}{2} \sqrt{\frac{|\alpha - \alpha'|}{\kappa(\zeta_0)}} \; e^{ik|\alpha - \alpha'|d(\mathbf{l})} \tag{8.21}$$

can be used for an approximate recovery of the curvature and the support function (modulo $2\pi / k|\alpha - \alpha'|$) of the obstacle, provided one knows that the total field satisfies the Dirichlet boundary condition. The uncertainty in the support function determination can be remedied by using different combinations of vectors $\alpha$ and $\alpha'$ as described in the numerical results section.

## 9. The Support Function Method for Neumann and Robin boundary conditions

Since it is also of interest to localize the obstacle in the case when the boundary condition is not a priori known, one can modify the SFM as follows.

For the Robin boundary condition

$$A(\alpha', \alpha) = \frac{e^{i\frac{\pi}{4}}}{\sqrt{8\pi k}}$$
$$\int_\Gamma \left\{ u(\mathbf{y}) \frac{\partial e^{-ik\alpha'\cdot\mathbf{y}}}{\partial N(\mathbf{y})} + hu(\mathbf{y})e^{-ik\alpha'\cdot\mathbf{y}} \right\} ds(\mathbf{y}) ,$$
(9.1)

and on $\Gamma_+$

$$u(\mathbf{x}) = \left(1 + \frac{ikN(\mathbf{x})\cdot\alpha + h}{ikN(\mathbf{x})\cdot\alpha - h}\right) e^{ik\mathbf{x}\cdot\alpha} .$$

In the Kirchhoff approximation one lets $u = 0$, $u_N = 0$ on $\Gamma_-$. Thus

$$A(\alpha', \alpha) = \frac{e^{i\frac{\pi}{4}}}{\sqrt{8\pi k}} \int_{\Gamma_+} \left(1 + \frac{ikN(\mathbf{y})\cdot\alpha + h}{ikN(\mathbf{y})\cdot\alpha - h}\right)$$
$$(-ikN(\mathbf{y})\cdot\alpha' + h)e^{ik(\alpha-\alpha')\cdot\mathbf{y}} ds(\mathbf{y})$$
$$= \frac{2ike^{i\frac{\pi}{4}}}{\sqrt{8\pi k}} \int_{\Gamma_+} \left(\frac{-ikN(\mathbf{y})\cdot\alpha' + h}{ikN(\mathbf{y})\cdot\alpha - h}\right) \quad (9.2)$$
$$N(\mathbf{y})\cdot\alpha \; e^{ik(\alpha-\alpha')\cdot\mathbf{y}} ds(\mathbf{y})$$

Now the Stationary Phase method applied to (9.2) gives

$$A(\alpha', \alpha) = \frac{1}{2}\sqrt{\frac{|\alpha-\alpha'|}{\kappa(\zeta_0)}} \; e^{ik|\alpha-\alpha'|d(\mathbf{l})} \left[1 + O(\frac{1}{k})\right]$$
(9.3)

as $k \to \infty$. Here $\kappa$ is the curvature of the boundary at the specular point.

Since in the resonance region the frequency $k$ is not large, one can argue that, under the assumption of a small curvature $\kappa$, the function

$$\mathbf{y} \to \frac{-ikN(\mathbf{y})\cdot\alpha' + h}{ikN(\mathbf{y})\cdot\alpha - h}$$

is a slowly changing one, and it can be approximated by the constant

$$C = \frac{-ikN(\mathbf{s}_0)\cdot\alpha' + h}{ikN(\mathbf{s}_0)\cdot\alpha - h} = \frac{ikN(\mathbf{s}_0)\cdot\alpha + h}{ikN(\mathbf{s}_0)\cdot\alpha - h}$$
$$= \frac{-ik\mathbf{l}\cdot\alpha + h}{-ik\mathbf{l}\cdot\alpha - h}$$

This is equivalent to replacing the total field $u(\mathbf{x})$ on $\Gamma_+$ by

$$u(\mathbf{x}) = \left(1 + \frac{ikN(\mathbf{s}_0)\cdot\alpha + h}{ikN(\mathbf{s}_0)\cdot\alpha - h}\right) e^{ik\mathbf{x}\cdot\alpha}$$

Therefore

$$A(\alpha', \alpha)$$
$$= \frac{2ikCe^{i\frac{\pi}{4}}}{\sqrt{8\pi k}} \int_{\Gamma_+} N(\mathbf{y})\cdot\alpha \; e^{ik(\alpha-\alpha')\cdot\mathbf{y}} ds(\mathbf{y})$$
(9.4)

Now the approximation by the Stationary Phase method gives

$$A(\alpha', \alpha) \sim \frac{C}{2}\sqrt{\frac{|\alpha-\alpha'|}{\kappa(\zeta_0)}} e^{ik(\alpha-\alpha')\cdot\mathbf{y}(\zeta_0)} \quad (9.5)$$

Since

$$C = \frac{ik|\alpha-\alpha'|-2h}{ik|\alpha-\alpha'|+2h} = e^{i(\pi-2\gamma)}$$

where

$$\gamma = \arctan\frac{k|\alpha-\alpha'|}{2h}$$

it follows that

$$A(\alpha', \alpha) \sim \frac{1}{2}\sqrt{\frac{|\alpha-\alpha'|}{\kappa(\zeta_0)}} \; e^{i(k|\alpha-\alpha'|d(\mathbf{l})-2\gamma+\pi)}$$

Let $\mathbf{l} \in S^1$ be fixed. Let

$$R(\mathbf{l}) = \{\alpha \in S^1 : |\alpha\cdot\mathbf{l}| > 1/\sqrt{2}\} .$$

Therefore $\sqrt{2} < |\alpha-\alpha'| \le 2$. In this range we approximate $\gamma$ by

$$\gamma_0 = \arctan \frac{k}{h}$$

and get

$$A(\alpha', \ \alpha) \sim \frac{1}{2}\sqrt{\frac{|\alpha - \alpha'|}{\kappa(\zeta_0)}} \ e^{i(k|\alpha - \alpha'|d(\mathbf{l}) - 2\gamma_0 + \pi)}$$

(9.6)

Now one can recover the Support Function $d(\mathbf{l})$ from (9.6), and the location of the obstacle.

## 10. Numerical results for the Support Function Method

In the first numerical experiment the obstacle is the circle

$$D = \{(x_1, x_2) \in \mathbb{R}^2 : (x_1 - 6)^2 + (x_2 - 2)^2 = 1\}.$$

(10.1)

It is reconstructed using the Support Function Method for two frequencies in the resonance region: $k = 1.0$, and $k = 5.0$. Table 3 shows how well the approximation (8.21) is satisfied for various pairs of vectors $\alpha$ and $\alpha'$ all representing the same vector $\mathbf{l} = (1.0, \ 0.0)$ according to (8.16). The Table shows the ratios of the approximate Scattering Amplitude $A_a(\alpha', \ \alpha)$ defined as the right hand side of the equation (8.21) to the exact Scattering Amplitude $A(\alpha', \ \alpha)$. Note, that for a sphere of radius $a$, centered at $x_0 \in \mathbb{R}^2$, one has

$$A(\alpha', \ \alpha) = -\sqrt{\frac{2}{\pi k}} \ e^{-i\frac{\pi}{4}} e^{ik(\alpha - \alpha') \cdot \mathbf{x}_0}$$

$$\sum_{l=-\infty}^{\infty} \frac{J_l(ka)}{H_l^{(1)}(ka)} \ e^{il(\theta - \beta)} \ ,$$

(10.2)

where $\alpha' = \mathbf{x}/|\mathbf{x}| = e^{i\theta}$, and $\alpha = e^{i\beta}$. Vectors $\alpha$ and $\alpha'$ are defined by their polar angles shown in Table 3.

TABLE 3. Ratios of the approximate and exact Scattering Amplitudes $A_\alpha(\alpha', \alpha) / A(\alpha', \alpha)$ for $\mathbf{1} = (1.0, 0.0)$.

| $\alpha'$ | $\alpha$ | $k = 1.0$ | $k = 5.0$ |
|---|---|---|---|
| $\pi$ | 0 | 0.88473–0.17487i | 0.98859–0.05846i |
| $23\pi/24$ | $\pi/24$ | 0.88272–0.17696i | 0.98739–0.06006i |
| $22\pi/24$ | $2\pi/24$ | 0.87602–0.18422i | 0.98446–0.06459i |
| $21\pi/24$ | $3\pi/24$ | 0.86182–0.19927i | 0.97997–0.07432i |
| $20\pi/24$ | $4\pi/24$ | 0.83290–0.22411i | 0.96701–0.08873i |
| $19\pi/24$ | $5\pi/24$ | 0.77723–0.25410i | 0.95311–0.10321i |
| $18\pi/24$ | $6\pi/24$ | 0.68675–0.27130i | 0.92330–0.14195i |
| $17\pi/24$ | $7\pi/24$ | 0.57311–0.25360i | 0.86457–0.14959i |
| $16\pi/24$ | $8\pi/24$ | 0.46201–0.19894i | 0.81794–0.22900i |
| $15\pi/24$ | $9\pi/24$ | 0.36677–0.12600i | 0.61444–0.19014i |
| $14\pi/24$ | $10\pi/24$ | 0.28169–0.05449i | 0.57681–0.31075i |
| $13\pi/24$ | $11\pi/24$ | 0.19019+0.00075i | 0.14989–0.09479i |
| $12\pi/24$ | $12\pi/24$ | 0.00000+0.00000i | 0.00000+-0.00000i |

Table 3 shows that only vectors $\alpha$ close to the vector $\mathbf{l}$ are suitable for the Scattering Amplitude approximation. This shows the practical importance of the backscattering data. As mentioned at the end of Section 8, any single combination of vectors $\alpha$ and $\alpha'$ representing $\mathbf{l}$ is not sufficient to uniquely determine the Support Function $d(\mathbf{l})$ from (8.21) because of the phase uncertainty. However, one can remedy this by using more than one pair of vectors $\alpha$ and $\alpha'$ as follows.

Let $\mathbf{l} \in S^1$ be fixed. Let

$$R(\mathbf{l}) = \{\alpha \in S^1 : |\alpha \cdot \mathbf{l}| > 1/\sqrt{2}\} \ .$$

Define $\Psi : \mathbb{R} \to \mathbb{R}^+$ by

$$\Psi(t) = \left\| \frac{A(\alpha',\alpha)}{|A(\alpha',\alpha)|} + e^{ik|\alpha - \alpha'|t} \right\|^2_{L^2(R(\mathbf{l}))}$$

where $\alpha' = \alpha'(\alpha)$ is defined by $\mathbf{l}$ and $\alpha$ according to (8.16), and the integration is done over $\alpha \in R(\mathbf{l})$ .

If the approximation (8.21) were exact for any $\alpha \in R(\mathbf{l})$, then the value of $\Psi(d(\mathbf{l}))$ would be zero. This justifies the use of the minimizer $t_0 \in \mathbb{R}$ of the function $\Psi(t)$ as an approximate value of the Support Function $d(\mathbf{l})$. If the Support Function is known for sufficiently many directions $\mathbf{l} \in S^1$, the obstacle can be localized using (8.6) or (8.10). The results of such a localization for $k = 1.0$ together with the original obstacle $D$ is shown on Figure 2. For $k = 5.0$ the identified obstacle is not shown, since it is practically the same as $D$. The only a priori assumption on $D$ was that it was located inside the circle of radius 20 with the center in the origin. The Support Function was computed for 16 uniformly distributed in $S^1$ vectors $\mathbf{l}$. The program run takes about 80 seconds on a 333 MHz PC.

In another numerical experiment we used $k = 1.0$ and a kite-shaped obstacle. Its boundary is described by

$$\mathbf{r}(t) = (5.35 + \cos t + 0.65 \cos 2t, \ 2.0 + 1.5 \sin t),$$
$$0 \le t < 2\pi \ .$$

(10.3)

Numerical experiments using the boundary integral equation method (BIEM) for the direct scattering problem for this obstacle centered in the origin are described in [10], section 3.5. Again, the Dirichlet boundary conditions were assumed. We computed the scattering amplitude for 120 directions $\alpha$ using the MRC method (see section ) with about 25% performance improvement over the BIEM, see [28].

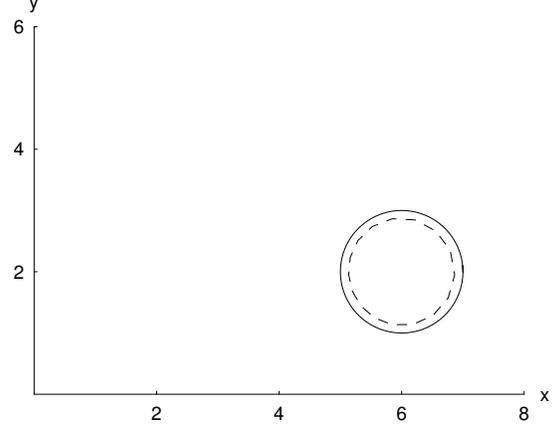

FIGURE 2. Identified (dotted line), and the original (solid line) obstacle $D$ for $k = 1.0$.

The Support Function Method (SFM) was used to identify the obstacle $D$ from the synthetic scattering amplitude with no noise added. The only a priori assumption on $D$ was that it was located inside the circle of radius 20 with the center in the origin. The Support Function was computed for 40 uniformly distributed in $S^1$ vectors $\mathbf{l}$ in about 10 seconds on a 333 MHz PC. The results of the identification are shown in Figure 3. The original obstacle is the solid line. The points were identified according to (8.10). As expected, the method recovers the convex part of the boundary $\Gamma$, and fails for the concave part. The same experiment but with $k = 5.0$ achieves a perfect identification of the convex part of the boundary. In each case the convex part of the obstacle was successfully localized. Further improvements in the obstacle localization using the MRC method are suggested in [27], and in the next section.

For the identification of obstacles with unknown boundary conditions let

$$A(t) = A(\alpha',\alpha) = |A(t)| e^{i\varphi(t)}$$

where, given $t$, the vectors $\alpha$ and $\alpha'$ are chosen as above, and the phase function

$\psi(t)$, $\sqrt{2} < t \leq 2$ is continuous. Similarly, let $A_a(t)$, $\psi_a(t)$ be the approximate scattering amplitude and its phase defined by formula (9.6).

If the approximation (9.6) were exact for an $\alpha \in R(\mathbf{l})$, then the value of $|\psi_a(t) - ktd(\mathbf{l}) + 2\gamma_0 - \pi|$ would be a multiple of $2\pi$.

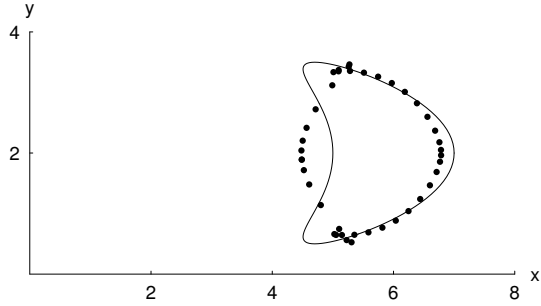

FIGURE 3. Identified points and the original obstacle $D$ (solid line); $k = 1.0$.

This justifies the following algorithm for the determination of the Support Function $d(\mathbf{l})$:

Use a linear regression to find the approximation $\psi(t) \approx C_1 t + C_2$ on the interval $\sqrt{2} < t \leq 2$. Then

$$d(\mathbf{l}) = \frac{C_1}{k} \quad . \qquad (10.4)$$

Also

$$h = -k \tan \frac{C_2}{2} \quad .$$

However, the formula for $h$ did not work well numerically. It could only determine if the boundary conditions were or were not of the Dirichlet type. Improvements are suggested in the next section.

TABLE 4. Identified values of the Support Function for the circle of radius 1.0 at $k = 3.0$.

| $h$ | Identified $d(\mathbf{l})$ | Actual $d(\mathbf{l})$ |
|---|---|---|
| 0.01 | –0.9006 | –1.00 |
| 0.10 | –0.9191 | –1.00 |
| 0.50 | –1.0072 | –1.00 |
| 1.00 | –1.0730 | –1.00 |
| 2.00 | –0.9305 | –1.00 |
| 5.00 | –1.3479 | –1.00 |
| 10.00 | –1.1693 | –1.00 |
| 100.00 | –1.0801 | –1.00 |

Table 4 shows that the algorithm based on (10.4) was successful in the identification of the circle of radius $1.0$ centered in the origin for various values of $h$ with no a priori assumptions on the boundary conditions. For this circle the Support Function $d(\mathbf{l}) = -1.0$ for any direction $\mathbf{l}$.

## 11. Inverse scattering methods based on the MRC

Suppose that an approximate location of the obstacle $D$ is obtained by a numerical inversion method, such as the SFM. Then one can try to use the MRC method to improve the localization of the boundary, see [27]. Such methods are under development by the authors, and they are going to be discussed elsewhere. Nevertheless, the MRC provides a tool for an easy construction of various examples illustrating the severe ill-posedness of the Inverse Scattering problem, which can be used for the algorithm's testing.

Here we present one such example. Let the obstacle $D$ be the unit circle $\{x \in \mathbb{R}^2 : |x| \leq 1\}$. If the incident field is $u_0(x) = e^{ikx\cdot\alpha}$, then the scattered field $v(x) = -u_0(x)$ for $x \in S = \partial D$, and its scattering amplitude is

$$A(\alpha', \alpha) = -\sqrt{\frac{2}{\pi k}} \, e^{-i\frac{\pi}{4}} \sum_{l=-\infty}^{\infty} \frac{J_l(ka)}{H_l^{(1)}(ka)} \, e^{il(\theta-\beta)} ,$$
(11.1)

where $\alpha' = \mathbf{x}/|\mathbf{x}| = e^{i\theta}$, and $\alpha = e^{i\beta}$.

Let $x_1 \in \mathbb{R}^2$. Fix an integer $L > 0$, and let $\mathbf{c} \in \mathbb{C}^{2L+1}$. Form the radiating solution

$$v_c(x) = \sum_{l=-L}^{L} c_l H_l^{(1)}(k|x-x_1|)e^{il\theta_1}, \quad (11.2)$$

where $(x-x_1)/|x-x_1| = e^{i\theta_1}$. Then its far field pattern is

$$A_{v_c}(\alpha') = \sqrt{\frac{2}{\pi k}}\ e^{-i\frac{\pi}{4}} \left( e^{-ik\alpha' \cdot x_1} \sum_{l=-L}^{L} c_l\ (-i)^l e^{il\theta} \right), \quad (11.3)$$

where $\alpha' = x/|x| = e^{i\theta}$.

TABLE 5. Near field values of two radiating solutions with practically the same far fields.

| $\alpha'$ | Re $v_c$ | Im $v_c$ | Re $v$ | Im $v$ |
|---|---|---|---|---|
| 0.00000 | –1189.60834 | –227.35213 | –0.54030 | –0.84147 |
| 0.31416 | –73.43878 | –15.81270 | –0.58082 | –0.81403 |
| 0.62832 | 1.94958 | 0.19051 | –0.69021 | –0.72361 |
| 0.94248 | 0.03298 | –0.52343 | –0.83217 | –0.55452 |
| 1.25664 | –1.07968 | –0.36021 | –0.95263 | –0.30412 |
| 1.57080 | –1.13445 | 0.00027 | –1.00000 | 0.00000 |
| 1.88496 | –0.96294 | 0.31629 | –0.95263 | 0.30412 |
| 2.19911 | –0.79021 | 0.55436 | –0.83217 | 0.55452 |
| 2.51327 | –0.66472 | 0.71819 | –0.69021 | 0.72361 |
| 2.82743 | –0.59154 | 0.81406 | –0.58082 | 0.81403 |
| 3.14159 | –0.56768 | 0.84565 | –0.54030 | 0.84147 |
| 3.45575 | –0.59154 | 0.81406 | –0.58082 | 0.81403 |
| 3.76991 | –0.66472 | 0.71819 | –0.69021 | 0.72361 |
| 4.08407 | –0.79021 | 0.55436 | –0.83217 | 0.55452 |
| 4.39823 | –0.96294 | 0.31629 | –0.95263 | 0.30412 |
| 4.71239 | –1.13445 | 0.00027 | –1.00000 | 0.00000 |
| 5.02655 | –1.07968 | –0.36021 | –0.95263 | –0.30412 |
| 5.34071 | 0.03298 | –0.52343 | –0.83217 | –0.55452 |
| 5.65487 | 1.94958 | 0.19051 | –0.69021 | –0.72361 |
| 5.96903 | –73.43878 | –15.8127 | –0.58082 | –0.81403 |

Fix an $\alpha \in S^1$, and let

$$r^{\min} = \min\{\|A_{v_c}(\alpha') - A(\alpha', \alpha)\| : \mathbf{c} \in \mathbb{C}^{2L+1}\}. \quad (11.4)$$

We conducted the minimization by the Singular Value Decomposition Method described in Section 3, using the following settings: $k = 1.0$, $L = 5$, $\alpha = (1.0, 0.0)$, and $x_1 = (0.8, 0.0)$. The $L^2$ norm in (11.4) was computed over $M = 120$ directions $\alpha'_m$ uniformly distributed in the unit circle $S^1$, and then normalized by $\sqrt{M}$, so that the identity function would have the norm equal to 1. The resulting value of the residual $r^{\min} = 0.00009776$ indicates that the far field $A(\alpha', \alpha)$ was practically perfectly fit by the radiating solution of the form (11.2). However, as the Table 5 shows, the restrictions of the exact scattered field $v$, and the fitted field $v_c$ to the boundary $S$ of the obstacle $D$ are vastly different. The columns in Table 5 correspond to the real and the imaginary parts of the scattered fields, and the rows correspond to different values of the angle $\alpha'$. Thus, one has to conclude that, as expected, a coincidence of the radiating solutions at the far field does not imply that the near fields are also coincidental.

## 12. Analysis of the Linear Sampling Method

During the last decade many papers were published, in which the obstacle identification methods were based on a numerical verification of the inclusion of some function $f := f(\alpha, z)$, $z \in \mathbb{R}^3$, $\alpha \in S^2$, in the range $R(B)$ of a certain operator $B$. Examples of such methods include [8], [9], [14]. However, one can show that the methods proposed in the above papers have essential difficulties, see [30]. Although it is true that $f \notin R(B)$ when $z \notin D$, it turns out that in any neighborhood of $f$ there are elements from $R(B)$. Also, although $f \in R(B)$ when $z \in D$, there are elements in every neighborhood of $f$ which do not belong to $R(B)$ even if $z \in D$. Therefore it is quite difficult to construct a stable numerical method for the identification of $D$ based on the verification of the inclusions $f \notin R(B)$, and $f \in R(B)$. Some published numerical results were intended to show that the method based on the above idea works practically, but it is not clear how these conclusions were obtained.

Let us introduce some *notations*: $N(B)$ and $R(B)$ are, respectively, the null-space and the range of a linear operator $B$, $D \in \mathbb{R}^3$ is a bounded domain (obstacle) with a smooth boundary $S$, $D' = \mathbb{R}^3 \setminus D$, $u_0 = e^{ik\alpha \cdot x}$, $k = const > 0$, $\alpha \in S^2$ is a unit vector, $N$ is the unit normal to $S$ pointing into $D'$, $g = g(x, y, k) := g(|x-y|) := \frac{e^{ik|x-y|}}{4\pi|x-y|}$, $f := e^{-ik\alpha' \cdot z}$, where $z \in \mathbb{R}^3$ and $\alpha' \in S^2$, $\alpha' := xr^{-1}$, $r = |x|$, $u = u(x, \alpha, k)$ is the scattering solution:

$$(\Delta + k^2)u = 0 \text{ in } D', \quad u|_S = 0, \quad (12.1)$$

$$u = u_0 + v, v = A(\alpha', \alpha, k)e^{ikr}r^{-1} + o(r^{-1}),$$
$$\text{as } r \to \infty, \quad (12.2)$$

where $A := A(\alpha', \alpha, k)$ is called the scattering amplitude, corresponding to the obstacle $D$ and the Dirichlet boundary condition. Let $G = G(x, y, k)$ be the resolvent kernel of the Dirichlet Laplacian in $D'$:

$$(\Delta + k^2)G = -\delta(x-y) \text{ in } D', \quad G|_S = 0, \quad (12.3)$$

and $G$ satisfies the outgoing radiation condition.

If
$$(\Delta + k^2)w = 0 \text{ in } D', \quad w|_S = h, \quad (12.4)$$
and $w$ satisfies the radiation condition, then ([23]) one has

$$w(x) = \int_S G_N(x, s)h(s)ds,$$
$$w = A(\alpha', k)e^{ikr}r^{-1} + o(r^{-1}), \quad (12.5)$$

as $r \to \infty$, and $xr^{-1} = \alpha'$. We write $A(\alpha')$ for $A(\alpha', k)$, and

$$A(\alpha') := Bh := \frac{1}{4\pi}\int_S u_N(s, -\alpha')h(s)ds, \quad (12.6)$$

as follows from Ramm's lemma:

**Lemma 1.** ([23], p.46) *One has:*

$$G(x, y, k) = g(r)u(y, -\alpha', k) + o(r^{-1}),$$
$$\text{as } r = |x| \to \infty, xr^{-1} = \alpha', \quad (12.7)$$

*where $u$ is the scattering solution of (12.1)-(12.2).*

One can write the scattering amplitude as:

$$A(\alpha', \alpha, k) = -\frac{1}{4\pi}\int_S u_N(s, -\alpha')e^{ik\alpha \cdot s}ds. \quad (12.8)$$

The following claim follows easily from the results in [23], [24] (cf [14]):

*Claim:* $f := e^{-ik\alpha' \cdot z} \in R(B)$ if and only if $z \in D$.

*Proof:* If $e^{-ik\alpha' \cdot z} = Bh$, then Lemma 1 and (12.6) imply

$$g(y, z) = \int_S G_N(s, y)hds \text{ for } |y| > |z|.$$

Thus $z \in D$, because otherwise one gets a contradiction: $\lim_{y \to z} g(y,z) = \infty$ if $z \in \overline{D'}$, while $\lim_{y \to z} \int_S G_N(s,y) h ds < \infty$ if $z \in \overline{D'}$. Conversely, if $z \in D$, then Green's formula yields $g(y,z) = \int_S G_N(s,y) g(s,z) ds$. Taking $|y| \to \infty, \frac{y}{|y|} = \alpha'$, and using Lemma 1, one gets $e^{-ik\alpha' \cdot z} = Bh$, where $h = g(s, z)$. The claim is proved.

Consider $B: L^2(S^2) \to L^2(S^2)$, and $A: L^2(S^2) \to L^2(S^2)$, where $B$ is defined in (12.6) and $Aq := \int_{S^2} A(\alpha', \alpha) q(\alpha) d\alpha$. Then one proves (see [30]):

**Theorem 1.** *The ranges $R(B)$ and $R(A)$ are dense in $L^2(S^2)$.*

**Remark 1.** In [8] a 2D inverse obstacle scattering problem is considered. It is proposed to solve the equation (1.9) in [8]:

$$\int_{S^1} A(\alpha, \beta) \varsigma d\beta = e^{-ik\alpha \cdot z}, \qquad (12.9)$$

where $A$ is the scattering amplitude at a fixed $k > 0$, $S^1$ is the unit circle, $\alpha \in S^1$, and $z$ is a point on $\mathbb{R}^2$. If $\varsigma = \varsigma(\beta, z)$ is found, the boundary $S$ of the obstacle is to be found by finding those $z$ for which $\|\varsigma\| := \|\varsigma(\beta, z)\|_{L^2(S^1)}$ is maximal. Assuming that $k^2$ is not a Dirichlet or Neumann eigenvalue of the Laplacian in $D$, and that $D$ is a smooth, bounded, simply connected domain, the authors state Theorem 2.1 [8], p. 386, which says that for every $\varepsilon > 0$ there exists a function $\varsigma \in L^2(S^1)$, such that

$$\lim_{z \to S} \|\varsigma(\beta, z)\| = \infty, \qquad (12.10)$$

and (see [8], p.386),

$$\|\int_{S^1} A(\alpha, \beta) \varsigma d\beta - e^{-ik\alpha \cdot z} \| < \varepsilon. \qquad (12.11)$$

There are several questions concerning the proposed method.

First, equation (12.9), in general, is not solvable. The authors propose to solve it approximately, by a regularization method. The regularization method applies for stable solution of solvable ill-posed equations (with exact or noisy data). If equation (12.9) is not solvable, it is not clear what numerical "solution" one seeks by a regularization method.

Secondly, since the kernel of the integral operator in (12.9) is smooth, one can always find, for any $z \in \mathbb{R}^2$, infinitely many $\varsigma$ with arbitrary large $\|\varsigma\|$, such that (12.11) holds. Therefore it is not clear how and why, using (12.10), one can find $S$ numerically by the proposed method.

A numerical implementation of the Linear Sampling Method (LSM) suggested in [8] consists of solving a discretized version of (12.9)

$$F\mathbf{g} = \mathbf{f}, \qquad (12.12)$$

where $F = \{A\alpha_i, \beta_j\}$, $i = 1,..., N$, $j = 1, 2, ... N$ is a square matrix formed by the measurements of the scattering amplitude for $N$ incoming, and $N$ outgoing directions. In 2-D the vector $\mathbf{f}$ is formed by

$$\mathbf{f}_n = \frac{e^{i\frac{\pi}{4}}}{\sqrt{8\pi k}} e^{-ik\alpha_n \cdot z}, \quad n = 1,..., N,$$

see [6] for details.

Denote the Singular Value Decomposition of the far field operator by $F = USV^H$. Let $s_n$ be the singular values of $F$, $\rho = U^H \mathbf{f}$, and $\mu = V^H \mathbf{f}$. Then the norm of the sought function $g$ is given by

$$\|\varsigma\|^2 = \sum_{n=1}^{N} \frac{|\rho_n|^2}{s_n^2} \qquad (12.13)$$

A different LSM is suggested by A. Kirsch in [14]. In it one solves

$$(F^*F)^{1/4} \mathbf{g} = \mathbf{f} \qquad (12.14)$$

instead of (12.12). The corresponding expression for the norm of $\varsigma$ is

$$\|\varsigma\|^2 = \sum_{n=1}^{N} \frac{|\mu_n|^2}{s_n} \qquad (12.15)$$

A detailed numerical comparison of the two LSMs and the linearized tomographic inverse scattering is given in [6].

The conclusions of [6], as well as of our own numerical experiments are that the method of Kirsch (12.14) gives a better, but a comparable identification, than (12.12). The identification is significantly deteriorating if the scattering amplitude is available only for a limited aperture, or the data are corrupted by noise. Also, the points with the *smallest* value of $\|\varsigma\|$ are the best in locating the inclusion, and not the *largest* one, as required by the theory in [14] and in [8]. In Figures 4 and 5 the implementation of the Colton-Kirsch LSM (12.13) is denoted by *gnck*, and of the Kirsch method (12.15) by *gnk*. The Figures show a contour plot of the logarithm of $\|\varsigma\|$. In all the cases the original obstacle was the circle of radius 1.0 centered at the point (10.0, 15.0). A similar circular obstacle that was identified by the Support Function Method (SFM) is discussed in Section 10. Note that the actual radius of the circle is 1.0, but it cannot be seen from the LSM identification. The LSM does not require any knowledge of the boundary conditions on the obstacle. The use of the SFM for unknown boundary conditions is discussed in Section 10. The LSM identification was performed for the scattering amplitude of the circle computed analytically with no noise added. In all the experiments the value for the parameter $N$ was chosen to be 128.

## 13. Remarks about other numerical methods for solving IOSP

There are numerous papers of various authors in which parameter-fitting procedures are proposed for solving IOSP. Some consider objective functionals that do not have an absolute minimum equal to zero, and even are not defined on the exact solution. For example, in [10] one takes a closed surface $S_0$ inside the unknown obstacle with the boundary $S$, assumes that $k^2$ is not a Dirichlet eigenvalue of the Laplacian in the domain $D_0$ with the boundary $S_0$, and uses the approximate relation

$$-4\pi A(\alpha', \alpha) = \int_{S_0} e^{-ik(\alpha', s)} h(s, \alpha) ds \qquad (13.1)$$

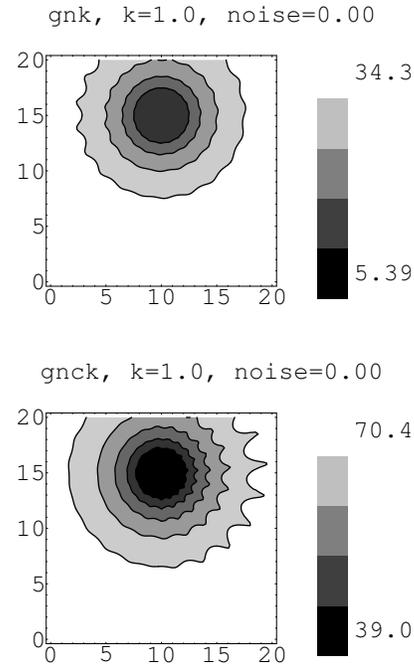

**FIGURE 4.** Identification of a circle at $k = 1.0$.

as an integral equation for the equation of $S$ in spherical coordinates $r = f(\alpha')$ (assuming $S$ to be star-shaped). One can prove that, for any fixed $\alpha \in S^2$, the set $\int_{S_0} e^{-ik(\alpha', s)} h(s, \alpha) ds$ is dense in $L^2(S^2)$ when $h$ runs through a dense subset of $L^2(S_0)$ (see [44]). In (13.1) the unknown is $h$. If $h$ is found, then the function $f$, which determines $S$, is to be calculated as the solution to the equation

$$\int_{S_0} g(f(\alpha')\alpha', s) h ds = -e^{ik(\alpha, f(\alpha')\alpha')} \quad (13.2)$$

where $g(x,s) := \dfrac{e^{ik|x-s|}}{4\pi|x-s|}$. It is recommended to solve equations (13.1) and (13.2) by a regularization method. In fact, equation (13.1) is, in general, not solvable, and equation (13.2) may be not solvable either.

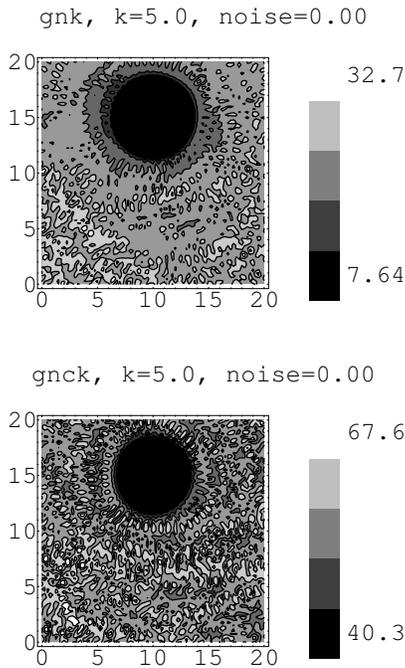

FIGURE 5. Identification of a circle at $k$=5.0.

The disadvantages of these methods are discussed in [25]. In [25] and [37] some objective functionals are proposed which are defined on the exact solutions of IOSP and attain absolute minimum equal to zero at the exact solutions.

Another approach is to use a regularized Newton-type method for solving nonlinear integral equation to which the IOSP can be reduced (e.g, see [15]). In [38] an inversion formula for calculating the obstacle from fixed-frequency scattering data is given, but an open problem is to find an algorithm for computing the function $\nu$ in this inversion formula. Such an algorithm for computing a similar function $\nu$ in the inverse potential scattering problem (IPSP) is given in [26].

## References

1. Alber, H. D., A quasi-periodic boundary value problem for the Laplacian and the continuation of its resolvent. Proc. Roy. Soc. Edinburgh Sect. A 82 (1978/79), no. 3-4, 251--272.
2. Albertsen N.C., Chesneaux J.-M., Christiansen S., Wirgin A., Comparison of four software packages applied to a scattering problem, Mathematics and Computers in Simulation, 48, (1999), 307-317.
3. Barantsev R, [1971] Concerning the Rayleigh hypothesis in the problem of scattering from finite bodies of arbitrary shapes, Vestnik Leningrad Univ., Math., Mech., Astron., **7**, 56-62.
4. A. Bonnet-Bendhia, Guided waves by eletromagnetic gratings and non-uniqueness examples for the diffraction problem, Math. Math. in the Appl. Sci., 17, (1994), 305-338.
5. A. Bonnet-Bendhia, K. Ramdani, Diffraction by an acoustic grating perturbed by a bounded obstacle, Adv. Comp. Math., 16, (2002), 113-138.
6. Brandfass M., Lanterman A.D., Warnick K.F., [2001] Acomparison of the Colton-Kirsch inverse scattering methods with linearized tomographic inverse scattering, Inverse Problems, **17**, 1797-1816.
7. S. Christiansen and R.E. Kleinman, On a misconception involving point collocation and the Rayleigh hypothesis, IEEE Trans. Anten. Prop., 44,10, 1309-1316, 1996. 760.
8. Colton, D., Kirsch, A., A simple method for solving inverse scattering problems in the resonance region, Inverse Problems 12 (1996), no. 4, 383--393
9. Colton, D., Coyle, J., Monk, P., Recent

developments in inverse acoustic scattering theory, SIAM Rev. 42 (2000), no. 3, 369--414
10. Colton D., Kress R. [1992] Inverse Acoustic and Electromagnetic Scattering Theory, Springer-Verlag, New York.
11. Eidus, D. M. Some boundary-value problems in infinite regions, Izv. Akad. Nauk SSSR Ser. Mat. 27, (1963) 1055--1080.
12. Kazandjian L, Rayleigh-Fourier and extinction theorem methods applied to scattering and transmission at a rough solid-solid interface, J. Acoust. Soc. Am., 92, 1679-1691, 1992.
13. Kazandjian L, Comments on "Reflection from a corrugated surface revisited", [J. Acoust. Soc. Am., 96, 1116-1129 (1994)]" J. Acoust. Soc. Am., 98, 1813-1814, (1995). 1245
14. Kirsch, A., Characterization of the shape of a scattering obstacle using the spectral data for far field operator, Inverse Probl., 14, (1998), 1489-1512
15. Kress, R., Integral equation methods in inverse obstacle scattering, ANZIAM J. 42 (2000), no. 1, 65--78.
16. Lions J.L., Magenes E. [1972] Non-Homogeneous Boundary Value Problems and Applications, Springer, New York.
17. Millar R. [1973] The Rayleigh hypothesis and a related least-squares solution to the scattering problems for periodic surfaces and other scatterers, Radio Sci., **8**, 785-796.
18. R. Millar, The Rayleigh hypothesis and a related least-squares solution to the scattering problems for periodic surfaces and other scatterers, Radio Sci., 8 (1973) 785-796.
19. R. Millar, On the Raleigh assumption in scattering by a periodic surface, Proc. Camb. Phil. Soc., 69 (1971) 217-225.; 65 (1969) 773-791.
20. Nazarov S., Plamenevskii B., Elliptic problems in domains with piecewise smooth boundaries. de Gruyter Expositions in Mathematics, 13. Walter de Gruyter, Berlin, 1994.
21. Petit R. (editor), Electromagnetic theory of gratings, Topics in Current Physics, 22. Springer-Verlag, Berlin-New York, 1980.
22. Press W.H., Teukolsky S.A., Vetterling W.T., Flannery B.P. [1992] Numerical Recepies in FORTRAN, Second Ed., Cambridge University Press.
23. Ramm A.G. [1986] Scattering by Obstacles, D. Reidel Publishing, Dordrecht, Holland.
24. Ramm A.G. [1992] Multidimensional Inverse Scattering Problems, Longman/Wiley, New York.
25. Ramm A.G. [1994] Multidimensional Inverse Scattering Problems, Mir, Moscow (expanded Russian edition of [24]).
26. Ramm A.G. [2002] Stability of solutions to inverse scattering problems with fixed-energy data, Milan Journ. of Math., 70, 97-161.
27. Ramm A.G. [2002] Modified Rayleigh Conjecture and Applications, J. Phys. A: Math. Gen. **35**, L357-L361.
28. Gutman S., Ramm A.G. [2002] Numerical Implementation of the MRC Method for obstacle Scattering Problems, J. Phys. A: Math. Gen. **35**, 8065-8074.
29. Ramm A.G., Gutman S. Modified Rayleigh Conjecture for Scattering by Periodic Structures, submitted.
30. Ramm A.G., Gutman S. Analysis of a Method for Identification of Obstacles, preprint.
31. A.G. Ramm, G. Makrakis, Scattering by obstacles in acoustic waveguides, Spectral and scattering theory, in the book: editor A.G. RAMM, Plenum publishers, New York, 1998, pp.89-110.
32. A.G. Ramm, Singularities of the inverses of Fredholm operators, Proc. of Roy. Soc. Edinburgh, 102A, (1986), 117-121.
33. A.G. Ramm, Investigation of the scattering problem in some domains with infinite boundaries I, II, Vestnik 7, (1963), 45-66; 19, (1963), 67-76.
34. A.G. Ramm, Uniqueness theorems for inverse obstacle scattering problems in Lipschitz domains, Applic. Analysis, 59, (1995), 377-383.
35. A.G. Ramm, Continuous dependence of the scattering amplitude on the surface of an obstacle, Math. Methods in the Appl. Sci.,

18, (1995), 121-126.
36. A.G. Ramm, Scattering amplitude as a function of the obstacle, Appl.Math.Lett., 6, N5, (1993), 85-87.
37. A.G. Ramm, Numerical method for solving inverse scattering problems, Doklady of Russian Acad. of Sci., 337, N1, (1994), 20-22
38. A.G. Ramm, Stability of the solution to inverse obstacle scattering problem, J. Inverse and Ill-Posed Problems, 2, N3, (1994), 269-275.
39. A.G. Ramm, Stability estimates for obstacle scattering, J. Math. Anal. Appl. 188, N3, (1994), 743-751.
40. A.G. Ramm, P.Pang and G.Yan, A uniqueness result for the inverse transmission problem, Int. Jour. of Appl. Math., 2, N5, (2000), 625-634.
41. A.G. Ramm, Reconstruction of the domain shape from the scattering amplitude, Radiotech. i Electron., 11, (1965), 2068-2070.
42. Gutman S., Ramm A.G. [2003] Support Function Method for Inverse Scattering Problems, In the book "Acoustics, mechanics and related topics of mathematical analysis", (ed. A.Wirgin), World Scientific, New Jersey, pp. 178-184.
43. A.G. Ramm, M.Sammartino, Existence and uniqueness of the scattering solutions in the exterior of rough domains, in the book "Operator Theory and Its Applications", Amer. Math. Soc., Fields Institute Communications vol.25, pp.457-472, Providence, RI, 2000. (editors A.G.Ramm, P.N.Shivakumar, A.V.Strauss).
44. Ramm, A.G., On a property of the set of radiation patterns, J. Math. Anal. Appl. 98, (1984), 92-98.
45. Gol'dshtein V., Ramm A.G., Embedding operators for rough domains, Math. Ineq. and Applic., 4, N1, (2001), 127-141.
46. Rayleigh J.W., On the dynamical theory of gratings, Proc. Roy. Soc. A, 79, (1907), 399-416.